\documentclass[12pt]{article}
\usepackage{amsmath}
\usepackage{latexsym}
\usepackage{amssymb}

\newcommand{\at}{\symbol{'100}}
\newcommand{\N}{{\mathbb N}}
\newcommand{\R}{{\mathbb R}}
\newcommand{\sub}{\subseteq}
\newcommand{\cg}{{\mathfrak g}}
\newcommand{\ch}{{\mathfrak h}}
\newcommand{\cF}{{\mathcal F}}
\newcommand{\mto}{\mapsto}

\newcommand{\dl}{{\displaystyle\lim_{\longrightarrow}}}
\DeclareMathOperator{\Diff}{Diff}

\DeclareMathOperator{\evol}{evol}
\DeclareMathOperator{\id}{id}
\begin{document}
$\;$\\[-27mm]
\begin{center}
{\Large\bf Direct limits of regular Lie groups}\\[3.9mm]
{\bf Helge Gl\"{o}ckner}
\end{center}
\begin{abstract}
\hspace*{-7.6mm}
Let $G$ be a regular Lie group which is a directed union of regular Lie groups $G_i$
(all modelled on possibly infinite-dimensional, locally convex spaces).
We show that $G=\dl\,G_i$\vspace{-.3mm}
as a regular Lie group if~$G$ admits a so-called direct limit chart.
Notably, this allows the regular Lie group $\Diff_c(M)$ of compactly supported diffeomorphisms
to be interpreted as a direct limit of the regular Lie groups $\Diff_K(M)$
of diffeomorphisms supported in compact sets $K\sub M$,
even if the finite-dimensional smooth manifold~$M$ is merely paracompact
(but not necessarily $\sigma$-compact), which is new.
Similar results are obtained for the
test function groups $C^k_c(M,F)$ with values in a Lie group~$F$. 
\end{abstract}
Consider a Lie group~$G$ (modelled on a locally convex space)
which is a union $G=\bigcup_{i\in I}G_i$ of such Lie groups
$G_i$, where $(I,\leq)$ is a directed set, $G_i\sub G_j$ for all $i\leq j$ in $I$
and all inclusion maps $\eta_{ji}\colon G_i\to G_j$ and $\eta_i\colon G_i\to G$
are smooth group homomorphisms.
It is natural to ask whether $G$ is the direct limit of the Lie groups $G_i$
in the category of Lie groups and smooth group homomorphisms;
or, equivalently, whether a group homomorphism
\[
f\colon G\to H
\]
to a Lie group~$H$ is smooth if $f|_{G_i}\colon G_i\to H$ is smooth for all $i\in I$.
This question and related ones (concerning direct limit properties of~$G$
as a topological group, as a smooth manifold, or as topological space)
are fairly well-understood if $I=\N$, so that we are dealing with a union
$G=\bigcup_{n\in\N}G_n$ of an ascending
sequence $G_1\sub G_2\sub\cdots$ of Lie groups (see \cite{CMP}
and the references therein).
For typical examples, consider
a paracompact finite-dimensional smooth manifold~$M$;
then the group $\Diff_c(M)$
of all smooth diffeomorphisms $\psi\colon M\to M$ such that $\psi(x)=x$
off some compact set is a Lie group.
Moreover, $\Diff_K(M):=\{\psi\in \Diff_c(M)\colon (\forall x\in M\setminus K)\;
\psi(x)=x\}$ is a Lie subgroup of $\Diff_c(M)$ for each compact subset $K\sub M$
(cf.\ \cite{Mic} and \cite{GaN}).
If $F$ is a Lie group with neutral element $e$ and $k\in \N_0\cup\{\infty\}$,
then the ``test tunction group'' $C^k_c(M,F)$ is a Lie group,
comprising all $C^k$-maps $\gamma\colon M\to F$ such that $\gamma(x)=e$
for $x\in M$ off some compact set~$K$ \cite{JFA,cplte,SEM,Mea}. For fixed $K$,
\[
C^k_K(M,F):=\{\gamma\in C^k_c(M,F)\colon (\forall x\in M\setminus K)\;\, \gamma(x)=e\}
\]
is a Lie subgroup of $C^k_c(M,F)$.
If $M$ with $\mbox{dim}(M)>0$ is $\sigma$-compact, non-compact and $K_1\sub K_2\sub \cdots$
is an exhaustion of $M$ by compact sets, then
\[
\Diff_{K_1}(M)\sub \Diff_{K_2}(M)\sub\cdots
\]
is an ascending sequence of Lie groups and
\[
\Diff_c(M)=\dl\, \Diff_{K_n}(M)\vspace{-1mm}
\]
holds as a Lie group and as a topological group, but not as a smooth manifold
(see \cite{CMP}), nor as a topological space (see \cite{TSH}). Likewise,
\[
C^\infty_c(M,F)
=\dl\, C^\infty_{K_n}(M,F)\vspace{-1mm}
\]
holds as a Lie group and topological group but neither as a smooth manifold
nor as a topological space if $M$ and the $K_n$ are as before
and $F$ is a non-discrete finite-dimensional Lie group.\\[2.3mm]
For ascending unions $G=\bigcup_{i\in I} G_i$ with uncountable index sets~$I$,
hardly anything is known concerning direct limit properties: neither general results,
nor results concerning concrete examples.
Notably, in the case of a paracompact finite-dimensional smooth
manifold~$M$ which fails to be $\sigma$-compact,
it is an open problem whether
\[
\Diff_c(M)=\dl\, \Diff_K(M)\vspace{-1mm}
\]
holds as a Lie group,
and whether
\[
C^k_c(M,F)=\dl\, C^k_K(M,F)\vspace{-1mm}
\]
holds as a Lie group
for each finite-dimensional Lie group~$F$ (see Problem~17.13
in the extended preprint version {\tt arxiv.math/0606078}
of \cite{CMP}).\\[2.3mm]
In this note, we explain that the situation improves
if we restrict attention to the class of Lie groups which are \emph{regular}
(in the sense recalled in Section~\ref{secprel}). Regularity is a key concept in infinite-dimensional
Lie theory going back to Milnor~\cite{Mil}
(in the case of sequentially complete modelling spaces);
see \cite{SEM}, \cite{GaN}, and \cite{Nee} for further information
(cf.\ also \cite{KaM2}). Up to now, no examples of non-regular Lie groups
modelled on sequentially complete (or Mackey complete)
locally convex spaces have been found.
Therefore, a focus on regular Lie groups
hardly poses a restriction.\\[2.3mm]
We consider a Lie group $G$ modelled on a locally convex spaces
which is the union $G=\bigcup_{i\in I}G_i$
of a directed family $(G_i)_{i\in I}$ of such Lie groups,
as described above.
We assume that $G$ ``has a direct limit chart'',
i.e., $G$ admits a chart around its neutral
element $e$ which is built up from compatible charts of the
Lie groups $G_i$ (see Section~\ref{secprel} for details). Then the following holds:\\[4mm]
{\bf Theorem~A.}
\emph{If the Lie group $G$ is an ascending union $G=\bigcup_{i\in I}G_i$
of Lie groups~$G_i$ and $G$ admits a direct limit chart,
then a group homomorphism $f\colon G\to H$ to a regular Lie group~$H$
is smooth if and only if $f|_{G_i}\colon G_i\to H$ is smooth for each $i\in I$.}\\[4mm]
We record an immediate consequence:\\[4mm]
{\bf Theorem~B.}
\emph{If a regular Lie group $G$ is an ascending union $G=\bigcup_{i\in I}G_i$
of regular Lie groups~$G_i$ and $G$ admits a direct limit chart,
then $G=\dl\, G_i$\vspace{-1mm} in the category of regular Lie groups
and smooth group homomorphisms.} $\,\square$\\[4mm]
Now consider a family $(G_j)_{j\in J}$ of Lie groups~$G_j$, where $G_j$ is modelled on the locally
convex space~$E_j$. Then the subgroup
\[
G:=\bigoplus_{j\in J}G_j:=\Big\{(g_j)_{j\in J}\in\prod_{j\in J} G_j\colon
\mbox{$g_j=e$ for all but finitely many $j$}\,\Big\}
\]
of $\prod_{j\in J}G_j$ can be made a Lie group in a natural way,
modelled on the locally convex direct sum $\bigoplus_{j\in J} E_j$ (see \cite{Mea},
where the notation $\prod_{j\in J}^*G_j$ is used for $\bigoplus_{j\in J} G_j$).
The Lie group $\bigoplus_{j\in J} G_j$ is called the \emph{weak direct product}
of the family $(G_j)_{j\in J}$.
If $\phi_j\colon V_j\to W_j\sub E_j$ is a chart of~$G_j$ with $e\in V_j$ and $\phi(e)=0$,
then $\bigoplus_{j\in J}V_j:=G\cap \prod_{j\in J}V_j$
is an open identity neighbourhood in~$G$ and the map
\begin{equation}\label{chawea}
\phi:=\oplus_{j\in J}\phi_j\colon\bigoplus_{j\in J} V_j\to\bigoplus_{j\in J} W_j\sub\bigoplus_{j\in J} E_j,
\;\, (g_j)_{j\in J}\mto (\phi_j(g_j))_{j\in J}
\end{equation}
is a chart for~$G$ (cf.\ \cite{Mea}).
Among other things, weak direct products are useful tools
for the study of diffeomorphism groups and test functions groups
(cf.\ \cite{CMP}, \cite{cplte}, \cite{Mea}, and \cite{GaN}).\\[4mm]
{\bf Theorem C (Direct limit properties of prime examples).}
\begin{itemize}
\item[(a)]
\emph{If $G:=\bigoplus_{j\in J} G_j$ is the weak direct product of
a family $(G_j)_{j\in J}$ of Lie groups,
then a group homomorphism $f\colon G\to H$ to a regular Lie group~$H$
is smooth if and only if $f|_{G_j}\colon G_j\to H$ is smooth for each $j\in J$.
Let ${\mathcal F}$ be the set of all finite subsets of~$J$ $($directed via inclusion$)$.
If $G$ and each $G_j$ is regular, then $G=\dl_{\Phi\in\cF}\prod_{j\in \Phi}G_j$\vspace{-.7mm}
holds in the category of regular Lie groups and smooth group homomorphisms.}
\end{itemize}
\emph{Now let $M$ be a paracompact, finite-dimensional smooth manifold.
Let ${\mathcal K}$ be the set of compact subsets of~$M$, directed via inclusion.
Then we have}:
\begin{itemize}
\item[(b)]
\emph{$\Diff_c(M)=\dl_{K\in {\mathcal K}}\Diff_K(M)$\vspace{-.4mm}
holds in the category of regular Lie groups and smooth group homomorphisms.}
\item[(c)]
\emph{If $F$ is a Lie group and $k\in \N_0\cup\{\infty\}$, then a group
homomorphism $f\colon C^k_c(M,F)\to H$ to a regular Lie group~$H$ is
smooth if and only if $f|_{C^k_c(M,F)}\colon C^k_K(M,F)\to H$ is
smooth for each $K\in{\mathcal K}$.
If $C^k_c(M,F)$ and each $C^k_K(M,F)$ is regular,
then $C^k_c(M,F)=\dl\, C^k_K(M,F)$\vspace{-.7mm}
in the category of regular Lie groups and smooth group homomorphisms.}\vspace{2mm}
\end{itemize}
Criteria ensuring that $\bigoplus_{j\in J}G_j$ is regular are known
(see \cite{SEM}), and recalled in Section~\ref{secprfC}. Notably, each weak direct
product of finite-dimensional Lie groups (or Banach-Lie groups) is regular.
Criteria ensuring that $C^k_c(M,F)$ and its Lie subgroups $C^k_K(M,F)$
are regular are known as well (see \cite{SEM})
and recalled in Section~\ref{secprfC}. Notably, $C^k_c(M,F)$ and $C^k_K(M,F)$
are regular whenever $F$ is a finite-dimensional Lie group
or a Banach-Lie group.\\[3mm]
\noindent
Let $G$ be a Lie group with tangent space $\cg:=L(G):=T_e(G)$ at the neutral element.
We say that $G$ has a \emph{smooth exponential function}
if there exists a (necessarily unique) smooth function $\exp_G\colon \cg\to G$
such that $\exp_G((t+s)v)=\exp_G(tv)\exp_G(sv)$ for all $v\in \cg$ and $s,t\in\R$,
and $\frac{d}{dt}\big|_{t=0}\exp_G(tv)=v$ for all $v\in \cg$.
If $G$ has a smooth exponential function and $\exp_G$ is a local $C^\infty$-diffeomorphism at~$0$,
then $G$ is called \emph{locally exponential}.
Every regular Lie group has a smooth exponential function.
Many important examples of Lie groups are locally exponential
(for example, all Lie groups modelled on Banach spaces),
but there also are many important examples which fail to be locally exponential
(like the diffeomorphism group $\Diff({\mathbb S})$ of the circle),
cf.\ \cite{Mil}, \cite{GaN}, \cite{Nee}. We record an easy observation:\\[4mm]
{\bf Theorem D.} \emph{Let $G$ be a locally exponential Lie group which is the ascending
union $G=\bigcup_{i\in I} G_i$ of Lie groups $G_i$ with smooth exponential function.
Assume that $L(G)=\dl\, L(G_i)$\vspace{-1mm} as a locally convex space
with limit maps $L(\eta_i):=T_e(\eta_i)$, for the direct system $((L(G_i)_{i\in I},(L(\eta_{ji})_{i\leq j})$
of locally convex spaces.
Then $G=\dl\,G_i$\vspace{-1mm} in the category of Lie groups with smooth exponential function
and smooth group homomorphisms.}\\[4mm]
Theorem~D was the starting point for this work. It is a general strategy of infinite-dimensional
Lie theory that one-parameter groups $t\mto \exp_G(tv)$ should be replaced with general
smooth curves $\gamma\colon [0,1]\to G$ with $\gamma(0)=e$ in classical proofs,
which are represented on the Lie algebra level
not by a single element $v\in \cg$, but by an element in $C^\infty([0,1], \cg)$
(the left logarithmic derivative of~$\gamma$).
Following this strategy, Theorems A and~B can be considered
as adaptations of Theorem~D when $G$ fails to be locally exponential,
so that smoothness of $f\circ\exp_G$ cannot guarantee smoothness of $f\colon G\to H$.\\[4mm]
{\bf Acknowledgement.}
The research was supported by
Deutsche Forschungsgemeinschaft (DFG), project GL 357/9-1.
\section{Preliminaries and Notation}\label{secprel}
Write $\N:=\{1,2,\ldots,\}$ and $\N_0:=\N\cup\{0\}$.\\[3mm]
{\bf 1.1} We shall work with $C^k$-maps between open subsets of locally convex
topological vector spaces as introduced by Bastiani~\cite{Bas} for $k\in \N_0\cup\{\infty\}$
(a setting which is also known as Keller's $C^k_c$-theory).
More generally, consider locally convex spaces $E$ and $F$
and a subset $U\sub E$ which is locally convex in the sense
that each point of~$U$ has a convex neighbourhood in~$U$.
Following~\cite{GaN}, we say that a map $f\colon U\to F$ is $C^k$
if $f$ is continuous and there exist continuous maps
\[
d^{(j)}f\colon U\times E^j\to F
\]
for all $j\in \N$ with $j\leq k$ such that the following
iterated directional derivatives exist and are given by the right hand side:
\[
(D_{y_j}\cdots D_{y_1}f)(x)=d^{(j)}f(x,y_1,\ldots,y_j)
\]
for all $x$ in the interior of~$U$ and $y_1,\ldots, y_j\in E$.
As usual, $C^\infty$-maps are also called \emph{smooth}.
We refer to \cite{RES} and \cite{GaN} for outlines of the
corresponding concepts of manifolds and Lie groups
modelled on arbitrary locally convex spaces
(which need not satisfiy any completeness conditions).
All Lie groups and manifolds we consider may be infinite-dimensional,
unless the contrary is stated.
Compare \cite{Mil} for the case of sequentially complete modelling spaces,
\cite{Ham} for differential calculus on Fr\'{e}chet spaces;
see also~\cite{Mic}.
If $U$ is an open subset of a locally convex space~$E$,
we identify its tangent bundle $TU$ with $U\times E$, as usual.
If $M$ is a smooth manifold, we let $\pi_{TM}\colon TM\to M$
be the bundle projection. If $f\colon M\to U$ is a $C^1$-map, we write
$df\colon TM\to E$ for the second component of $Tf\colon TM\to U\times E$;
thus $Tf=(f\circ\pi_{TM},df)$.\\[3mm]
{\bf 1.2} Every Lie group $G$ acts smoothly on its tangent bundle from the left via
\[
G\times TG\to TG,\quad (g,v)\mto g.v:=T\lambda_g(v)\quad\mbox{for $g\in G$ and $v\in TG$,}
\]
where $\lambda_g\colon G\to G$, $x\mto gx$ is left translation by~$g$.
Write $\cg :=L(G):=T_e(G)$ for the tangent space of~$G$ at the neutral element.
If $\gamma\colon [0,1]\to G$ is a $C^1$-curve, we define its left logarithmic derivative as
\[
\delta(\gamma)\colon [0,1]\to\cg,\quad t\mto\gamma(t)^{-1}.\gamma'(t),
\]
where $\gamma'(t):=T\gamma(t,1)$. Endow $C([0,1],\cg)$
with the compact-open topology. Given $k\in \N_0\cup\{\infty\}$,
we equip $C^k([0,1],\cg)$ with the so-called \emph{compact-open $C^k$-topology}
(the initial topology with respect to the maps $C^k([0,1],\cg)\to C([0,1],\cg)$,
$\gamma\mto \gamma^{(j)}$ taking $\gamma$ to its $j$th derivative
for all $j\in \N_0$ with $j\leq k$).
If each $\gamma\in C^k([0,1],\cg)$ arises as $\gamma=\delta(\eta)$
for a (necessarily unique) $C^{k+1}$-map $\eta\colon [0,1]\to G$
with $\eta(0)=e$ and the map
\[
\mbox{evol}\colon C^k([0,1], \cg)\to G,\quad \gamma\mto\eta(1)
\]
is smooth, then $G$ is called \emph{$C^k$-regular}.
We mention that $C^k$-regularity implies $C^\ell$-regularity for all $\ell\in\N\cup\{\infty\}$
such that $\ell\geq k$.
If $G$ is $C^\infty$-regular (the weakest regularity property),
we simply say that~$G$ is \emph{regular}.
See \cite{Mil}, \cite{SEM}, \cite{GaN}, and \cite{Nee} for further information
(cf.\ also \cite{KaM2}).\\[3mm]
{\bf 1.3} Consider a Lie group~$G$ which is the union of a directed system $(G_i)_{i\in I}$
of Lie groups, as described in the introduction, with inclusion maps $\eta_i\colon G_i\to G$
for $i\in I$ and $\eta_{ji}\colon G_i\to G_j$ for $i\leq j$ in~$I$.
Let $E$ be the modelling space of~$G$ and $E_i$ be the modelling space of~$G_i$,
for $i\in I$.
We say that a chart $\phi\colon V\to W\sub E$ of the smooth manifold~$G$
is a \emph{direct limit chart} if $e\in V$ and the following holds:
\begin{itemize}
\item[(a)]
$E=\bigcup_{i\in I} E_i$ and
$E_i\sub E_j$ for $ i\leq j$ in~$I$; moreover the inclusion maps $\lambda_i\colon E_i\to E$
and $\lambda_{ji}\colon E_i\to E_j$ are continuous linear and
$E$, with the maps $\lambda_i$, is the direct limit locally convex space
$E=\dl\, E_i$\vspace{-.7mm}
of the direct system $((E_i)_{i\in I},(\lambda_{ji})_{i\leq j})$ of locally convex spaces.
\item[(b)]
There are open $e$-neighbourhoods $V_i\sub G_i$ and open $0$-neighbourhoods
$W_i\sub E_i$ such that $\phi(V_i)=W_i$ and $\phi_i:=\phi|_{V_i}\colon V_i\to W_i$
is a chart for~$G_i$, for each $i\in I$, and moreover
$V=\bigcup_{i\in I} V_i$
and $V_i\sub V_j$ (and hence also $W_i\sub W_j)$ for all $i\leq j$ in~$I$.
\end{itemize}
Then the $L(\eta_i):=T_e(\eta_i)$ and $L(\eta_{ji}):=T_e(\eta_{ji})$ are injective and $L(G)$,
with maps $L(\eta_i)$, is the locally convex direct limit of
$((L(G_i))_{i\in I},(L(\eta_{ji}))_{i\leq j})$.
When convenient, we identify $L(G)$ with $E$ and $L(G_i)$ with $E_i$
using the isomorphisms $d\phi|_{L(G)}$ and $d\phi_i|_{L(G_i)}$.\\[2.3mm]
See \cite{SUR} for generalities concerning direct limits of Lie groups
and related topics. In the case of ascending sequences $G_1\sub G_2\sub\cdots$,
direct limit charts were defined and exploited in~\cite{CMP}.
The case of uncountable index sets was considered in~\cite{Hom}.
\section{Proof of Theorem~A}
To prove Theorem~A, we
write $\cg_i:=L(G_i)$ for $i\in I$ and
abbreviate $\cg:=L(G)$ and $\ch:=L(H)$.
For $i\in I$, let $\eta_i\colon G_i\to G$ be the inclusion map
and for $i\leq j$ in $I$,
let $\eta_{ji}\colon G_i\to G_j$ be the inclusion map; all of these are group homomorphisms and smooth.
Abbreviate $\psi_i:=L(f|_{G_i})\colon \cg_i\to \ch$.
Since $f|_{G_i}=f|_{G_j}\circ\eta_{ji}$, we have $\psi_i=\psi_j\circ L(\eta_{ji})$
for all $i\leq j$ in~$I$. As $\cg=\dl\, \cg_i$ as a locally convex space,
we deduce that there is a unique continuous linear map
$\psi\colon \cg\to \ch$ such that
\[
\psi\circ L(\eta_i)=\psi_i\quad\mbox{for all $\, i\in I$.}
\]
Let $\phi\colon V\to W\sub \cg$ be a direct limit chart for~$G$.
Thus $V=\bigcup_{i\in I} V_i$ with open $e$-neighbourhoods $V_i\sub G_i$
and $W=\bigcup_{i\in I}W_i$ with open $0$-neighbourhoods $W_i\sub\cg_i$
such that $W_i\sub W_j$ and $V_i\sub V_j$ for all $i\leq j$ in~$I$
(identifying $\cg_i$ with $L(\eta_i)(\cg_i)\sub \cg$ as a vector space)
and $\phi|_{V_i}\colon V_i\to W_i$ is a chart for~$G_i$.
We let $Q\sub W$ be an open $0$-neighbourhood which is balanced
(i.e., $[{-1},1] Q\sub Q$). After replacing $W$ with~$Q$, the set $V$ with $\phi^{-1}(Q)$,
the set $W_i$ with $W_i\cap Q$ and $V_i$ with $\phi^{-1}(W_i\cap Q)$,
we may assume that~$W$ is balanced.
Then
\begin{equation}\label{lineWi}
(\forall w\in W)(\exists j\in I)\;\, [0,1]w\sub W_j.
\end{equation}
In fact, let $w\in W$. For each $t\in [0,1]$, we have $tw\in W$ and
find $i(t)\in I$ such that $tw\in W_{i(t)}$. As the map $[0,1]\to \cg_{i(t)}$, $s\mto sw$
is continuous and $W_{i(t)}\sub \cg_{i(t)}$ is open, $t$ has a neighbourhood $U(t)\sub [0,1]$
such that $U(t)w\sub W_{i(t)}$.
By compactness of $[0,1]$, there exist $t_1,\ldots, t_n\in [0,1]$ such that
$[0,1]=U(t_1)\cup\cdots\cup U(t_n)$. Let $j\in I$ such that $i(t_k)\leq j$
for all $k\in\{1,\ldots, n\}$. Then $U(t_k)w\sub W_{i(t_k)}\sub W_j$ for each $k$
and thus $[0,1]w\sub W_j$ (establishing (\ref{lineWi})).\\[2.3mm]
Given $x\in V$, we define a smooth curve $\gamma_x\colon [0,1]\to G$ with $\gamma_x(0)=e$ and
$\gamma_x(1)=x$ via
\[
\gamma_x\colon [0,1]\to G,\quad t\mto \phi^{-1}(t\phi(x)).
\]
Now $[0,1]\phi(x)\sub W_j$ for some $j\in I$, by (\ref{lineWi}).
Then $\gamma_x(t)=(\phi|_{V_j}^{W_j})^{-1}(t\phi(x))\in G_j$
for all $t\in [0,1]$, and this is a smooth $G_j$-valued function of~$t$.
Thus
\[
(\forall x\in V)(\exists j\in I)\; \gamma_x([0,1])\sub G_j\;\,\mbox{and}\;\,
\gamma_x|^{G_j}\in C^\infty([0,1],G_j).
\]
For $x$ and $j$ as before, we obtain
\begin{eqnarray}
f(x)&=& (f\circ\gamma_x)(1)=(f|_{G_j}\circ \gamma_x|^{G_j})(1)
=\evol_H(\delta(f|_{G_j}\circ \gamma_x|^{G_j}))\notag\\
&=& \evol_H(\psi_j\circ \delta(\gamma_x|^{G_j}))
=\evol_H(\psi\circ L(\eta_j)\circ\delta(\gamma_x|^{G_j}))\notag\\
&=&\evol_H(\psi\circ \delta(\eta_j\circ \gamma_x|^{G_j}))
=\evol_H(\psi\circ\delta(\gamma_x)).\label{longcompo}
\end{eqnarray}
Note that $h\colon V\times [0,1]\to W\sub \cg$, $(x,t)\mto t\phi(x)$ is a smooth map, whence also
\[
h^\vee\colon V\to C^\infty([0,1],W)\sub C^\infty([0,1],\cg), \quad
x\mto h(x,\cdot)
\]
is smooth (see \cite[Theorem~B]{AaS}). Now $C^\infty([0,1],V)$ is an open identity neighbourhood
in $C^\infty([0,1],G)$ and the map
\[
C^\infty([0,1],\phi^{-1})\colon C^\infty([0,1],W)\to C^\infty([0,1],V),
\;\;
\gamma\mto \phi^{-1}\circ\gamma
\]
is a $C^\infty$-diffeomorphism by the construction of the Lie group structure on $C^\infty([0,1],G)$
(see, e.g., \cite{GaN}; cf.\ \cite{JFA}).
Now $\delta\colon C^\infty([0,1],G)\to C^\infty([0,1],\cg)$, $\gamma\mto\delta(\gamma)$
is smooth (see \cite[Lemma~2.1]{SEM})
and $C^\infty([0,1],\cg)\to C^\infty([0,1],\ch)$, $\gamma\mto \psi\circ\gamma$
is smooth (and continuous linear). Since $\phi^{-1}\circ h^\vee(x)=\gamma_x$ for each $x\in V$,
using (\ref{longcompo}) we see that
\[
f|_V=\evol_H\circ \,C^\infty([0,1],\psi)\circ \delta\circ C^\infty([0,1],\phi^{-1})\circ h^\vee
\]
is smooth. Like every homomorphism between Lie groups
which is smooth on an open identity neighbourhood, $f$ is smooth.\,$\square$
\section{Proof of Theorem C}\label{secprfC}
Proof of~(a). The map $\phi$ from (\ref{chawea}) is a direct limit chart for
$\bigoplus_{j\in J}G_j=\bigcup_{\Phi\in {\mathcal F}}\prod_{j\in \Phi}G_j$,
as $\bigoplus_{j\in J} V_j=\bigcup_{\Phi\in {\mathcal F}}\prod_{j\in\Phi}V_j$,
$\bigoplus_{j\in J}W_j=\bigcup_{\Phi\in {\mathcal F}}\prod_{j\in \Phi}W_j$
and $\phi$ restricts to the chart
\[
\prod_{j\in \Phi}\phi_j\colon\prod_{j\in \Phi}V_j\to\prod_{j\in\Phi}W_j
\]
of $\prod_{j\in\Phi}G_j$ around~$e$, for each $\Phi\in{\mathcal F}$.
Also note that the restriction of~$f$ to $\prod_{j\in\Phi}G_j$
with $\Phi=\{j_1,\ldots, j_n\}$
is the map
\[
\prod_{j\in\Phi}G_j\to H,\quad (g_j)_{j\in\Phi}\mto f|_{G_{j_1}}(g_{j_1})\cdot
f|_{G_{j_2}}(g_{j_2})\cdot\ldots\cdot  f|_{G_{j_n}}(g_{j_n}),\vspace{-.4mm}
\]
which is smooth if and only if $f_j$ is smooth for all $j\in \Phi$.
Therefore all assertions follow from Theorems A and B.

(b) Write $\pi_{TM}\colon TM\to M$ for the bundle projection
and let $\Gamma_c(TM)$
be the locally convex space of all compactly supported smooth vector fields on~$M$.
Given a compact set $K\sub M$, write $\Gamma_K(TM)$ for the
Fr\'{e}chet space of all smooth vector fields $X\colon M\to TM$
which are supported in~$K$.
Let $\Sigma\colon U\to M$ be a local addition for~$M$, i.e.,
a smooth map on an open neighbourhood $U\sub TM$ of $0_M:=\{0_p\in T_pM\colon p\in  M\}$
such that $\Sigma(0_p)=p$ for all $p\in M$ and moreover the map
\[
\theta\colon U\to M\times M, \quad v\mto (\pi_{TM}(v),\Sigma(v))
\]
has open image and is a $C^\infty$-diffeomorphism onto its image.\footnote{It is well-known
that such local additions always exist; one can take the Riemannian exponential map
for a Riemannian metric on~$M$ and restrict it to a suitable open set~$U$.}
We can (and shall) assume, moreover, that $T\Sigma|_{U\cap T_pM}(0,\cdot)=\id_{T_pM}$
for all $p\in M$.
There is an open subset
\[
\Omega\sub \{X\in \Gamma_c(TM)\colon X(M)\sub U\}
\]
containing $0_M$
such that $\Sigma\circ X\in \Diff_c(M)$ for all $X\in \Omega$ and the map
\[
\phi\colon \Omega\to \Diff_c(M),\quad X\mto \Sigma\circ X
\]
has open image $\Omega'$ and is a $C^\infty$-diffeomorphism onto its image.
Then
$\Omega_K:=\Omega\cap \Gamma_K(TM)$
is an open subset of $\Gamma_K(TM)$ and
the restriction $\phi_K$ to $\Omega_K$ has open image $\Omega_K'=\Omega'\cap\Diff_K(M)$
and is a $C^\infty$-diffeomorphism onto the latter. Identifying
$\Gamma_c(TM)$ with $T_{\id_M}(\Diff_c(M))$ by means of $T\phi(0,\cdot)$
and $\Gamma_K(TM)$ with $T_{\id_M}(\Diff_K(M))$ by means of $T\phi_K(0,\cdot)$,
the inclusion map $\Diff_K(M)\to \Diff_c(M)$ has the inclusion map
$\Gamma_K(TM)\to\Gamma_c(TM)$ as its tangent map at~$\id_M$.
Since
\[
\Gamma_c(TM)=\dl_{K\in{\mathcal F}} \Gamma_K(TM)
\]
as a locally convex space, we see that $\phi^{-1}$ is a direct limit chart for $\Diff_c(M)$.
As $\Diff_c(M)$ and the Lie groups $\Diff_K(M)$ are regular (see, e.g., \cite{Mea}), we can
apply Theorems~A and~B.

(c) Let $E$ be the modelling space of~$F$ and $\phi\colon V\to W\sub E$
be a chart for $F$ such that $e\in V$, $\phi(e)=0$ and $d\phi|_{L(F)}=\id_{L(F)}$.
Then $C^k_c(M,V)$ is an open identity-neighbourhood in $C^k_c(M,F)$;
moreover, $C^k_c(M,W)$ is an open $0$-neighbourhood in the locally convex space
$C^k_c(M,E)$ and
\[
C^k_c(M,\phi)\colon C^k_c(M,V)\to C^k_c(M,W),\quad\gamma\mto \phi\circ\gamma
\]
is a chart for $C^k_c(M,F)$. For each compact subset $K\sub M$, this chart restricts to the chart
\[
C^k_K(M,\phi)\colon C^k_K(M,V)\to C^k_K(M,W),\quad \gamma\mto\phi\circ \gamma
\]
of $C^k_K(M,F)$. We now identify $T_e(C^k_c(M,F))$ with $C^k_c(M,E)$
by means of the restriction of $dC^k_c(M,\phi)$ to an isomorphism between
the latter. Likewise, using $dC^k_K(M,\phi)$ we identify $T_e(C^k_K(M,F))$
with $C^k_K(M,E)$. Then the tangent map of the inclusion map $C^k_K(M,F)\to C^k_c(M,F)$
is the inclusion map $C^k_K(M,E)\to C^k_c(M,E)$. Since
\[
C^k_c(M,E)=\dl_{K\in{\mathcal K}} C^k_K(M,E)
\]
as a locally convex space, we see that $C^k_c(M,\phi)$ is a direct limit chart
for $C^k_c(M,F)$. Thus Theorem~A applies and if $C^k_c(M,F)$
and each of the Lie groups $C^k_K(M,F)$ is assumed regular, then also
Theorem~B applies. $\,\square$\\[4mm]
\noindent
{\bf Remark.} The regularity requirements in Theorem~C (a) and~(c)
are satisfied in the following situations:\medskip

(a) If $k\in \N_0$ and $G_j$ is $C^k$-regular for each $j\in J$,
then $\bigoplus_{j\in J}G_j$ is a $C^{k+1}$-regular Lie group
(see \cite[Corollary 13.6]{SEM})
and hence regular.\medskip\smallskip

(b) If $F$ is a regular Lie group, then $C^k_K(M,F)$ is
regular (this can be shown like \cite[Proposition~12.1]{SEM}).\medskip\smallskip

(c) If $F$ is $C^k$-regular for some $k\in \N_0$, then
$C^k_c(M,F)$ is $C^{k+1}$-regular (see \cite[Proposition~12.3]{SEM})
and hence regular.\\[2.3mm]
Classes of Lie groups which are $C^k$-regular for finite~$k$
can be found in \cite{SEM} and \cite{Mea}.\footnote{All Lie groups
which are measurably regular are $C^0$-regular.}
Notably, every Banach-Lie group is $C^0$-regular (and hence each finite-dimensional Lie group).
Moreover, all direct limits of ascending sequences of finite-dimensional
Lie groups are $C^0$-regular, and also all of the Lie groups $\Diff_c(M)$ and $\Diff_K(M)$.
\section{Proof of Theorem D}
Let $H$ be a Lie group with smooth exponential function and $f\colon G\to H$ be a homomorphism
of groups such that $f|_{G_i}\colon G_i\to H$ is smooth for each $i\in I$.
Then $L(f|_{G_i})=L(f|_{G_j}\circ\eta_{ji})=L(f|_{G_j})\circ L(\eta_{ji})$
for all $i\leq j$ in~$I$. By the universal property of the direct limit, there is a unique continuous linear map
$\psi\colon L(G)\to L(H)$ such that $\psi\circ L(\eta_i)=L(f|_{G_i})$ for all $i\in I$.
For each $i\in I$,
\[
\exp_H\!\circ \psi\circ L(\eta_i)=\exp_H\!\circ L(f|_{G_i})=f|_{G_i}\circ \exp_{G_i}
=f\circ\eta_i\circ\exp_{G_i}=f\circ\exp_G\circ L(\eta_i)
\]
and thus $\exp_H(\psi(v))=f(\exp_G(v))$ for all $v\in L(\eta_i)(L(G_i))$.
As $L(G)=\bigcup_{i\in I}L(\eta_i)(L(G_i))$, we deduce that $f\circ \exp_G=\exp_H\circ\psi$,
which is a smooth map. As $\exp_G$ is a local $C^\infty$-diffeomorphism at~$0$,
we deduce that~$f$ is smooth on some open identity neighbourhood in~$G$.
Since $f$ is a homomorphism, smoothness of~$f$ follows. $\,\square$.
{\small

{\small{\bf Helge  Gl\"{o}ckner}, Institut f\"{u}r Mathematik, Universit\"at Paderborn,\\
Warburger Str.\ 100, 33098 Paderborn, Germany.
Email: {\tt  glockner\at{}math.upb.de}}}\vfill
\end{document}